# GOODNESS-OF-FIT TESTING AND QUADRATIC FUNCTIONAL ESTIMATION FROM INDIRECT OBSERVATIONS

### By Cristina Butucea

### *Université Paris X and Université Paris VI*

We consider the convolution model where i.i.d. random variables $X_i$ having unknown density $f$ are observed with additive i.i.d. noise, independent of the $X$'s. We assume that the density $f$ belongs to either a Sobolev class or a class of supersmooth functions. The noise distribution is known and its characteristic function decays either polynomially or exponentially asymptotically.

We consider the problem of goodness-of-fit testing in the convolution model. We prove upper bounds for the risk of a test statistic derived from a kernel estimator of the quadratic functional $\int f^2$ based on indirect observations. When the unknown density is smoother enough than the noise density, we prove that this estimator is $n^{-1/2}$ consistent, asymptotically normal and efficient (for the variance we compute). Otherwise, we give nonparametric upper bounds for the risk of the same estimator.

We give an approach unifying the proof of nonparametric minimax lower bounds for both problems. We establish them for Sobolev densities and for supersmooth densities less smooth than exponential noise. In the two setups we obtain exact testing constants associated with the asymptotic minimax rates.

**1. Introduction.** We consider the *convolution model*,

$$(1) \qquad\qquad Y_i = X_i + \varepsilon_i, \qquad i = 1, \ldots, n,$$

where the random variables $X_i, \varepsilon_i$ are independent. We denote the common unknown density of $X_i$, $i = 1, \ldots, n$, by $f$. Let $\Phi(u) = \int e^{ixu} f(x)\,dx$ denote its characteristic function. We observe only the $Y_i$, $i = 1, \ldots, n$.

We consider the following nonparametric classes of density functions $f : \mathbb{R} \to \mathbb{R}_+$ with $\int f = 1$ and belonging to $\mathbb{L}_2$. A *Sobolev class* of density functions

---

with smoothness $\beta > 0$ and radius $L > 0$ is defined by

$$(2) \qquad W(\beta, L) = \left\{ f \in \mathcal{C}^\beta, \int |\Phi(u)|^2 |u|^{2\beta} \, du \le 2\pi L \right\}.$$

A class of *supersmooth* density functions for $\alpha, r, L > 0$, constants, is defined by

$$(3) \qquad S(\alpha, r, L) = \left\{ f \in \mathcal{C}^\infty, \int |\Phi(u)|^2 \exp(2\alpha |u|^r) \, du \le 2\pi L \right\}.$$

Let the noise be i.i.d. with known probability density $g$ and characteristic function $\Phi^g$. Then the resulting observations have common density $p = f \star g$ and characteristic function $\Phi^p = \Phi \cdot \Phi^g$. We also consider noise having a nonnull Fourier transform, $\Phi^g(u) \ne 0$, $\forall\, u \in \mathbb{R}$. Typically two different behaviors are distinguished in nonparametric estimation, *polynomially smooth* (or polynomial) noise

$$(4) \qquad |\Phi^g(u)| \sim |u|^{-\sigma}, \qquad |u| \to \infty,\ \sigma > 1,$$

and *exponentially smooth* (or supersmooth or exponential) noise

$$(5) \qquad |\Phi^g(u)| \sim \exp(-\gamma |u|^s), \qquad |u| \to \infty,\ \gamma, s > 0.$$

The first problem considered in this paper is nonparametric minimax goodness-of-fit testing from noisy data; that is, for a given density $f_0$ in the smoothness class $W(\beta, L_0)$, respectively, $S(\alpha, r, L_0)$ with $L_0 < L$, decide whether

$$H_0 : f = f_0,$$

or

$$H_1(\mathcal{C}, \psi_n) : f \text{ is in the smoothness class } \int (f - f_0)^2 \ge \mathcal{C} \psi_n^2,$$

from observations $Y_1, \ldots, Y_n$, for some fixed $\mathcal{C} > 0$ and $\psi_n > 0$.

Many important applications of this problem can be found in biology, medicine and physics, where errors-in-variables models have been extensively used.

In genomics, it is appropriate to admit that microarray data contain errors from non-biological sources. Gene expression is measured by scanning the fluorescence intensity of the microarray (see, e.g., Speed [29]). Software packages give slightly different results due to different correction and normalization methods. Testing the underlying fluorescence density from the scanned measurements provides a calibration method to the practitioner's particular microarray and scanner.

In medicine, many measurements are known to be subject to additive errors. In particular, the National Health and Nutrition Examination Survey (1976–1980), NHANES II, is a large dataset source of many studies of errors-in-variables models; see, for example, Carroll, Ruppert and Stefanski [7] for a previous study, NHANES I, and Delaigle and Gijbels [9]. The log-daily saturated fat intake is known to be a typical variable subject to error measurement and its probability density was estimated in the convolution model, with errors having either a Laplace or a Gaussian law. This variable is used to predict breast cancer, so the study is limited to women aged from 25 to 50. It was noted that the underlying density is symmetric, very smooth and has tails heavier than a normal distribution. Goodness-of-fit testing would help to choose between different types of densities.

Another important application of our testing procedure is to mixing location families $\{g(\cdot - \theta)\}_\theta$ with unknown mixing probability density $f$. The observation $Y$ therefore has probability density $p(y) = \int g(y - \theta) f(\theta) \, d\theta = f \star g$, as in the convolution model.

Moreover, we suggest use of this methodology for determining $K$, the unknown number of components in a finite mixture model. The astronomy dataset from Roeder and Wasserman [28], consisting of velocities ($\times 10^{-2}$) at which 82 galaxies from Corona Borealis spread away from our galaxy, was thoroughly studied in a $K$-mixture model with unknown $K$; see, for example, Stephens [31] and Richardson and Green [27]. Let $\theta_1, \ldots, \theta_K$ be the unknown states with the finite mixing probabilities $\{p_1, \ldots, p_K\}$. In order to fit into our theoretical framework, we suggest replacing the finite probability by a continuous law having density $f_{0K} = p_k \sum_{k=1}^K f_0(\cdot - \theta_k)$, with $f_0$ a peaked, supersmooth density. A preliminary estimation for different values of $K = 1, \ldots, \bar{K}$ provides estimators for $\{\theta_k\}_{k=1,\ldots,K}$ and $\{p_k\}_{k=1,\ldots,K}$. Then use goodness-of-fit testing as described later to test $H_0 : f = f_{0K}$ iteratively for $K = \bar{K}, \ldots, 1$ until the null is accepted.

All the previous examples, among many other applications, fit our setting for different values of parameters associated with the underlying density and the noise. These examples were treated from the point of view of estimating the deconvolution density, not that of the testing problem. To our knowledge this is the first time minimax testing is performed from data contaminated with errors. We give here simulation results showing very good testing properties between densities of the same families. As expected, testing quality is improved as the noise distribution becomes less smooth and/or has smaller variance. The test statistic has amazing convergence quality.

In the convolution model (1), the problem of nonparametric estimation of the deconvolution density $f$ has been intensively studied over the past two decades. In this paper, in order to surpass difficulties of estimation we address different issues, principally the goodness-of-fit test from noisy data in the $\mathbb{L}_2$ norm.

DEFINITION 1.   For a given $0 < \xi < 1$, a test statistic $\Delta_n^*$ is said to attain the testing rate $\psi_n$ over the smoothness class if there exists $\mathcal{C}^* > 0$ such that

$$(6) \qquad \limsup_{n \to \infty} \left\{ P_{f_0}[\Delta_n^* = 1] + \sup_{f \in H_1(\mathcal{C}, \psi_n)} P_f[\Delta_n^* = 0] \right\} \leq \xi$$

for all $\mathcal{C} > \mathcal{C}^*$. The rate $\psi_n$ is called the minimax rate of testing, if there exists $\mathcal{C}_* > 0$ and

$$(7) \qquad \liminf_{n \to \infty} \inf_{\Delta_n} \left\{ P_{f_0}[\Delta_n = 1] + \sup_{f \in H_1(\mathcal{C}, \psi_n)} P_f[\Delta_n = 0] \right\} \geq \xi$$

for all $0 < \mathcal{C} < \mathcal{C}_*$, where the inf is taken over all test procedures $\Delta_n$.

Moreover, if $\mathcal{C}^* = \mathcal{C}_*$ we call $\psi_n$ the exact (or sharp) minimax rate of testing.

We recall that the usual procedure is to construct the test statistic $\Delta_n^*$ such that (6) holds, also called the upper bound of the testing rate, and then prove the minimax optimality of this procedure, that is, the lower bounds in (7). If the test procedure does not depend on the smoothness of the unknown functions (which may vary in some interval), it is called adaptive to the smoothness and $\psi_n$ is the minimax adaptive rate.

Minimax and adaptive theory of testing has been extensively developed in density, regression and Gaussian white noise models when direct observations are available. For nonparametric minimax rates in goodness-of-fit testing in different setups we refer to Ingster [18], Ermakov [11] and references therein. Exact minimax rates have been found; see, for example, Lepski and Tsybakov [22] for the regression model with pointwise and sup-norm distances. The first adaptive rates were given by Spokoiny [30]. Exact minimax rates of testing for supersmooth functions are known only in the case $r = 1$ and for the Gaussian white noise model (see Pouet [26]) with pointwise and sup-norm distances. A further development consists of a goodness-of-fit test for a parametric composite null hypothesis and adaptive to the smoothness as in Fromont and Laurent [14] and Gayraud and Pouet [15]. Goodness-of-fit tests can be based on the distribution function rather than the density function of our data. In view of results by Fan [12] the $n^{-1/2}$ rates are still not feasible when estimating the distribution function in the convolution model. In view of numerous practical applications of testing, we expect the same problem in the context of data contaminated with errors to find similar extensive use in applied problems.

Here, the goodness-of-fit problem is considered in quadratic norm, $(\int (f - f_0)^2)^{1/2}$. As we can expect, the testing problem is easier than deconvolution density estimation, that is, the testing rates are faster as they appear in

Table 2. Note that minimax $\mathbb{L}_2$ testing can be performed at nearly the parametric rate $(\log n)^{(\sigma+1/4)/r} n^{-1/2}$ for supersmooth densities and polynomial noise.

We actually give exact minimax rates of testing in setups with densities less regular than the noise: Sobolev densities and exponential noise, supersmooth densities less smooth than the corresponding exponential noise ($r < s$).

The natural test statistic in this context is an estimator of $\int (f - f_0)^2$, where $f_0$ is given, from noisy data. Therefore, the second important problem treated in this paper is the estimation of the quadratic functional $d := \int f^2$, where $f$ is the density in the convolution model (1).

DEFINITION 2. An estimator $d_n$ of $d$ is said to attain the rate $\varphi_n$ over the smoothness class $W(\beta, L)$, respectively, $S(\alpha, r, L)$, if there exists a constant $C > 0$ such that

$$(8) \qquad \limsup_{n \to \infty} \sup_f \varphi_n^{-1} E_f[|d_n - d|] \le C,$$

and this rate is called minimax if no other estimator attains better rates uniformly over the class

$$(9) \qquad \liminf_{n \to \infty} \inf_{\hat{d}_n} \sup_f \varphi_n^{-1} E_f[|\hat{d}_n - d|] \ge c$$

for some $c > 0$, depending on fixed known parameters, where the supremum is taken over all densities in the smoothness class and the infimum over all estimators $\hat{d}_n$.

In some cases $n^{-1/2}$-consistent estimators of $d$ exist and we prove the asymptotic efficiency Cramér–Rao bound for such estimators (also called efficient estimators).

DEFINITION 3. An estimator $d_n$ of $d = \int f^2$ is asymptotically normally distributed with asymptotic variance $\mathcal{W} = \mathcal{W}(f)$ if

$$\sqrt{n}(d_n - d) \xrightarrow{d} N(0, \mathcal{W}(f)).$$

Moreover, it attains the asymptotic efficiency Cramér–Rao bound if for any $f_0$ in the Sobolev class $W(\beta, L)$, respectively, in $S(\alpha, r, L)$, and a family of shrinking neighborhoods $\mathcal{V}(f_0)$ of $f_0$,

$$\inf_{\mathcal{V}(f_0)} \liminf_{n \to \infty} \sup_{f \in \mathcal{V}(f_0)} n E_f[(\hat{d}_n - d)^2] \ge \mathcal{W}(f_0)$$

for any other estimator $\hat{d}_n$ of $d$.

6 C. BUTUCEA

When direct observations are available, it is well established that parametric rates can be achieved for smooth enough densities belonging, for example, to the Hölder class. Lower bounds for slower rates were found by Bickel and Ritov [1] for smoothness values less than $1/4$. In this context, Laurent [20] gave efficient estimation at the parametric rate and Birgé and Massart [2] proved nonparametric lower bounds for estimating more general quadratic functionals. The study of general functionals was completed by Kerkyacharian and Picard [19] for minimax rates and Tribouley [32] for adaptive estimation. Nemirovski [25] gave asymptotically efficient estimators of less smooth functionals, one or two times continuously differentiable.

In this paper, we give minimax results for setups in the nonparametric "regime" and efficiency constant in the sense of the theory of Ibragimov and Khas'minskii [17] and Levit [23] for asymptotically normal, $n^{-1/2}$-consistent estimators (see Table 1).

Moreover, it is possible to generalize these results to models with partially known noise distribution. Following results by Butucea and Matias [5], we can consider noise distributions with unknown scaling parameter (some more assumptions are needed in order to insure identifiability in the model). Current work is addressing the question of finding test procedures that will require even less information about the noise distribution.

These procedures can also be made adaptive, that is, free of the smoothness parameters, in some setups. We conjecture a loss of $\sqrt{\log n}$ due to adaptation to $\beta$ for estimating $d$ (see Efromovich and Low [10]), respectively, $\sqrt{\log \log n}$ for testing in the setup of Sobolev classes and polynomial noise. On the contrary, the testing procedure can be made fully data dependent with no loss in the rate for Sobolev densities and exponential noise and we expect the same to happen for estimating $d$. For supersmooth densities, computing the loss for adaptation is still an open problem.

The structure of the paper is as follows. In Section 2 we introduce the estimator $d_n$ of $\int f^2$ and the test statistic $\Delta_n^*$ and give some simulation results. In Section 3 we indicate the choice of the bandwidth in a functional's estimator in order to prove either upper bounds in the minimax sense, or asymptotic normality and efficiency, according to different setups. In Section 4 we deal with the goodness-of-fit testing problem and, for each setup, we compute upper bounds for testing rates. Finally, in Section 5 we describe the approach unifying the proofs of minimax nonparametric lower bounds from Sections 3 and 4 and prove them for nonparametric setups of Sobolev classes of densities and polynomial, respectively, exponential noise, and for the bias dominated setup of supersmooth densities less smooth than exponentially smooth noise ($r < s$). We have provided detailed proofs for one setup (Sobolev densities and polynomial noise) and put all other proofs in the Appendix that the interested reader may find in a longer version of this paper [4].

**2. Methodology and numerical results.**   In the described model, we consider the problem of estimating $d = \int f^2$, from available observations $(Y_i)_{i=1,\dots,n}$, where the density $f$ of observations $(X_i)_{i=1,\dots,n}$ is unknown. Let us denote the deconvolution kernel $K_n$ defined via its Fourier transform as

$$(10) \qquad \Phi^{K_n}(u) = \left( \Phi^g \left( \frac{u}{h} \right) \right)^{-1} \Phi^K(u),$$

where $K(x) = \sin(x)/(\pi x)$ is such that $\Phi^K(u) = I_{[|u| \le 1]}$ and the bandwidth $h = h_n \to 0$ when $n \to \infty$ will be specified later.

Define $d_n$, a bias-reduced estimator of $d$, by

$$(11) \qquad d_n = \frac{1}{n(n-1)} \sum_{k \ne j=1}^{n} \int K_{n,h}(x - Y_k) K_{n,h}(x - Y_j) \, dx,$$

where $K_{n,h}(\cdot) = 1/h K_n(\cdot/h)$.

In the sequel, we denote the $\mathbb{L}_2$ scalar product of two functions $M$ and $N$ by $\langle M, \overline{N} \rangle = \int M(x) N(x) \, dx$ and the complex conjugate of $N$ by $\overline{N}$.

In direct models, such a kernel based estimator can be found in Hall and Marron [16]. A biased-reduced kernel estimator first appeared in Bickel and Ritov [1], who proved that it is efficient for Hölder type smoothness values greater than $1/4$. Projection estimators were defined in Fan [13], Efromovich and Low [10] and Laurent [20].

Let us construct a test statistic from noisy data. It is natural to suggest as a test statistic $|T_n^*|$ the optimal estimator of the quadratic functional $\|f - f_0\|_2^2$,

$$T_n^* = \frac{1}{n(n-1)} \sum_{k \ne j} \langle K_{n,h}(\cdot - Y_k) - f_0, K_{n,h}(\cdot - Y_j) - f_0 \rangle,$$

where $h \searrow 0$ with $n$ and $K_n$ is defined in (10).

Define the test procedure

$$(12) \qquad\qquad \Delta_n^* = I[|T_n^*| > \mathcal{C}^* t_n^2]$$

for a constant $\mathcal{C}^* > 0$ and some threshold $t_n > 0$ depending on the setup.

In this paper, we chose the *sinc* kernel $K$, which has optimality properties. We stress the fact that for numerical implementation better choices are available, as was discussed in Butucea and Tsybakov [6]. Indeed, truncation of the Fourier transform gives a kernel $K_n$ which has $\int |K_n| = \infty$. It is enough to smooth $\Phi^K$ into a continuous trapezoidal-shaped function to get an absolutely integrable kernel. We actually use

$$\Phi^K(u) = I[|u| \le 1] + \exp(1 - (|u|(2 - |u|))^{-2}) \cdot I[1 \le |u| \le 2],$$

an infinitely differentiable function with compact support. The resulting deconvolution kernel has as many finite moments as $g$, the density of the noise, and the same optimality properties as our kernel $K_n$.

We consider $N = 100$ samples of size $n = 500$ and estimate the first type of error and the power of our testing procedure, as well as the mean squared error of our test statistic for estimating $\|f - f_0\|_2^2$.

The noise distribution will be either ordinary smooth $Laplace(1) \times S + M$ having density $g(x) = 0.5/S \exp(-|x - M|/(2S))$ and characteristic function $\Phi^g(u) = e^{iuM}(1 + (Su)^2)^{-1}$, of order $s = 2$, or $Laplace(3) \times S + M$ obtained as the sum of three independent rescaled $Laplace(1)$ variables, of order $s = 6$.

Densities $f_0$ under the null hypothesis are either Gaussian $N(M, S)$, belonging to a class of supersmooth functions $S(\alpha, r, L)$ with $r = 2$ and $\alpha < S^2/2$, or the Laplace density $Laplace(10) \times S + M$ having characteristic function $\Phi_0 = (1 + (Su)^2/(10))^{-10}$ belonging to a Sobolev class with $\beta < 20 - 1/2$. Other examples of densities can be found in Comte, Rozenholc and Taupin [8], including $Gamma$, $\chi^2$, stable distributions, densities with compactly supported characteristic functions and their mixtures.

We tested $f_0$: $N(1, 1)$ against successively rescaled Gaussian laws $N(1, 1 + (i - 1) \times 0.25)$, $i = 1, \ldots, 8$, under both $Laplace(1)$ and $Laplace(3)$ errors. We also tested $f_0$: $Laplace(10) \times 2$ against rescaled $Laplace(10) \times (2 + (i - 1) \times 0.25)$ and $f_0$: $N(1, 1)$ against shifted $N(1 + (i - 1) \times 0.25, 1)$, for $i = 1, \ldots, 8$, under $Laplace(1) \times \sqrt{0.5}$ errors.

We get excellent estimated test power, rapidly increasing with $i = 1, \ldots, 8$. We note that the power of the tests improves with the smoothness of tested densities, but it degrades with the smoothness of errors (when the signal to noise ratio is constant). These tests benefit from remarkable estimation properties of the test statistic $T_n^*$, as we can see from the boxplots in Figure 1.

We also note that the results are very satisfactory for detecting a one-mode density against a mixture of two identical densities. On the contrary, it is difficult to detect a heavier tailed density than $f_0$ when all other parameters are identical. This is due to the choice of the $\mathbb{L}_2$ norm, and this drawback is known in the testing literature. It would therefore be interesting and it is still an open problem to design tests with different distances ($\mathbb{L}_\infty$, Kullback or $\chi^2$ distance in the alternative) in this model.

**3. Estimation of $\int f^2$ in the convolution model.** In this section we present convergence properties of $d_n$ in (11) together with corresponding optimal choice of tuning parameters in each setup. Rates are summed up in Table 1.

DEFINITION 4. Let $d_n$ in (11) be the estimator of $d$ with bandwidth $h > 0$. We call the bias and the variance of this estimator, respectively,

$$B(d_n) \overset{\Delta}{=} |E_f[d_n] - d| \quad \text{and} \quad V(d_n) \overset{\Delta}{=} E_f[|d_n - E_f[d_n]|^2].$$

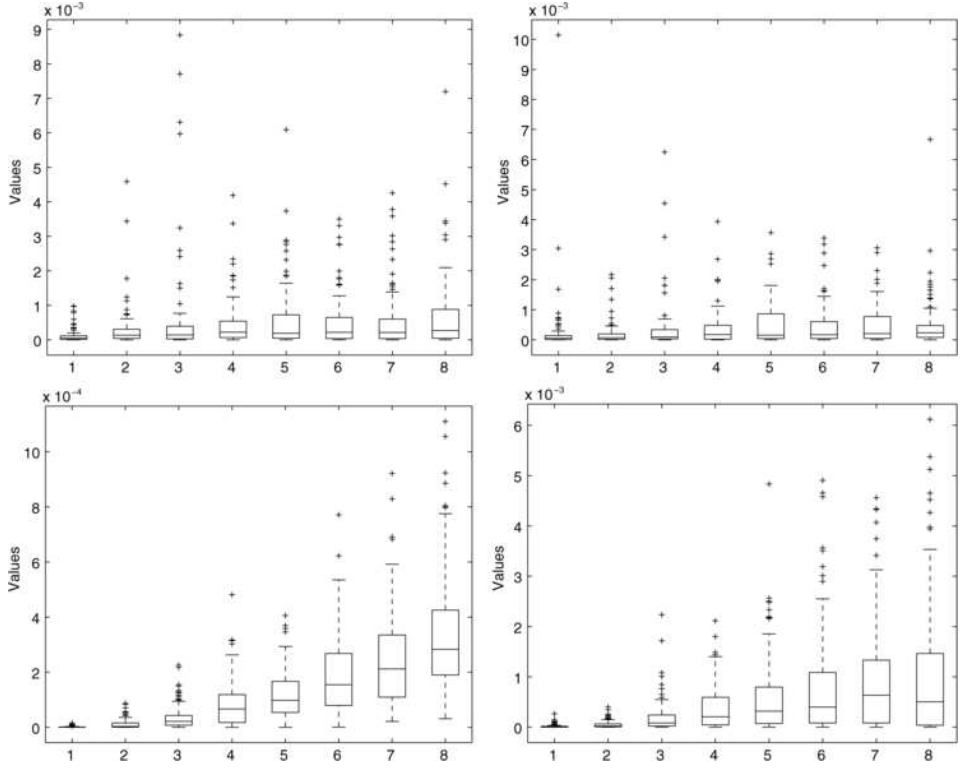

Fig. 1. *Mean square error of $T_n^*$ for estimating $\|f - f_0\|_2^2$, for $f_0$: $N(1,1)$ and $f$: $N(1, 1 + (i-1) \times 0.25)$ with Laplace(1) and Laplace(3) errors, in the upper graphics; for $f_0$: Laplace(10) $\times$ 2 and $f$: Laplace(10) $* (2 + (i-1) \times 0.5)$ and $f_0$: $N(1,1)$ and $f$: $N(1 + (i-1) \times 0.25, 1)$ with Laplace(1) $\times \sqrt{0.5}$ errors in the lower graphics, $i = 1, \ldots, 8$. The scaling is either $\times 10^{-3}$ or $\times 10^{-4}$, as indicated above the panels.*

**3.1. Sobolev densities and polynomial noise.** We study in detail the case where the underlying density $f$ belongs to a Sobolev class $W(\beta, L)$, with $\beta, L > 0$, defined in (2), and the noise is polynomial as defined in (4).

PROPOSITION 1. *For any density function $f$ in the Sobolev class $W(\beta, L)$, the estimator $d_n$ in (11) with bandwidth $h > 0$, $h \to 0$ as $n \to \infty$ is such that*

$$B(d_n) \leq L h^{2\beta},$$

$$V(d_n) = \frac{2\|p\|_2^2}{n^2 h^{4\sigma+1}} \frac{(1+e_n)}{\pi(4\sigma+1)} + \frac{4\Omega_g^2(f)}{n}(1+e_n)I_{\beta \geq \sigma} + \frac{E_n}{nh^{2(\sigma-\beta)+1}}I_{\beta < \sigma},$$

*where $\Omega_g(f) \geq 0$ is defined later in (13) and the sequences $e_n$ and $E_n$ do not depend on $f$ but depend only on $\beta$, $L$ and the noise density $g$, such that $e_n \to 0$ as $n \to \infty$, and $E_n$ bounded.*

In order to define $\Omega_g(f)$, let us note that for any $f$ in the Sobolev class $W(\beta, L)$ and $g$ a noise density satisfying (4) with $\beta \geq \sigma$, we have $\Phi/\overline{\Phi^g}$ a continuous function which is absolutely and quadratically integrable (see Lemma 4). Then we can define the function

$$F(y) = \frac{1}{2\pi} \int e^{-iyu} \frac{\Phi(u)}{\overline{\Phi^g}(u)} \, du,$$

which is real-valued and uniformly continuous, but not necessarily a density function. It is known (see Lukacs [24]) that if both characteristic functions $\Phi$ and $\overline{\Phi^g}$ are analytic around 0, then their quotient cannot be the characteristic function of any distribution function. Nevertheless, this function is bounded and its $L_2$ norm is uniformly bounded over densities $f$ in the Sobolev class by $M^F$ depending only on $\beta$, $L$ and the fixed given density $g$.

Let

$$
\begin{aligned}
\Omega_g^2(f) &\overset{\text{def}}{=} \int F^2(y) p(y) \, dy - \left( \int f^2(x) \, dx \right)^2 \\
&= E_f[F^2(Y)] - (E_f[F(Y)])^2.
\end{aligned}
\tag{13}
$$

Indeed, $\int f(x)^2 \, dx = (2\pi)^{-1} \langle \Phi, \overline{\Phi} \rangle = (2\pi)^{-1} \langle \Phi^p, \overline{\Phi}/\Phi^g \rangle = \langle p, F \rangle = E_f[F(Y)]$, which is therefore a real number.

REMARK 1. Note that (13) says that $4\Omega_g^2(f) = 4V(F(Y))$. This is heuristically similar to the results by given Laurent [20]. She estimates $\int f^2$ from direct observations and obtains the efficiency constant $4V(f(X)) = 4 \int f^3 - 4(\int f^2)^2$ when $\beta \geq 1/4$. In Theorems 1 and 2 we describe the same change of "regime" when $\beta \geq \sigma + 1/4$, respectively, $\beta < \sigma + 1/4$. Similarities between deconvolution with $\sigma$-polynomial noise and the derivative of order $\sigma$ have been noticed before. Indeed, we actually estimate $\int f^2 = E_f[F(Y)]$ here, where $F \star g = f$, whenever the function $F$ exists, and $F$ is as difficult to estimate as the $\sigma$-derivative of the function $f$.

PROOF OF PROPOSITION 1. Let us note that

$$
\begin{aligned}
E_f[d_n] &= E_f[\langle K_{n,h}(\cdot - Y_1), K_{n,h}(\cdot - Y_2) \rangle] \\
&= \|K_{n,h} \star p\|_2^2 = \|K_h \star f\|_2^2 \\
&= \frac{1}{2\pi} \int \Phi^K(hu) |\Phi(u)|^2 \, du.
\end{aligned}
\tag{14}
$$

By the Plancherel formula and equation (14)

$$
\begin{aligned}
B(d_n) &= \frac{1}{2\pi} \left| \int (\Phi^K(hu) - 1) |\Phi(u)|^2 \, du \right| \\
&\leq \frac{1}{2\pi} \int_{|u| > 1/h} (h|u|)^{2\beta} |\Phi(u)|^2 \, du \leq L h^{2\beta}.
\end{aligned}
$$

As for the variance let us first write

$$d_n - E_f[d_n] = \frac{1}{n(n-1)} \sum_{k \neq j}^{n} \langle K_{n,h}(\cdot - Y_k) - K_h \star f, K_{n,h}(\cdot - Y_j) - K_h \star f \rangle$$

$$+ \frac{2}{n} \sum_{k=1}^{n} \langle K_{n,h}(\cdot - Y_k) - K_h \star f, K_h \star f \rangle = S_1 + S_2, \quad \text{say}.$$

The variables in $S_1$ and in $S_2$ are uncorrelated and all of them are centered. Thus, $V(d_n) = E_f[S_1^2] + E_f[S_2^2]$. We have

$$E_f[S_1^2] = \frac{2}{n(n-1)} (E_f[\langle K_{n,h}(\cdot - Y_1) - K_h \star f, K_{n,h}(\cdot - Y_2) - K_h \star f \rangle^2])$$

$$= \frac{2}{n(n-1)} (E_f[\langle K_{n,h}(\cdot - Y_1), K_{n,h}(\cdot - Y_2) \rangle^2]$$

$$- 2E_f[\langle K_{n,h}(\cdot - Y_1), K_h \star f \rangle^2] + \|K_h \star f\|_2^4).$$

Moreover, $E_f[S_2^2] = 4n^{-1}(E_f[\langle K_{n,h}(\cdot - Y_1), K_h \star f \rangle^2] - \|K_h \star f\|_2^4)$. Similarly to Butucea [3], we have

$$E_f[\langle K_{n,h}(\cdot - Y_1), K_{n,h}(\cdot - Y_2) \rangle^2]$$

$$= \frac{1}{h} \int \int \frac{1}{h} \left| \int K_n\left(z + \frac{v-u}{h}\right) K_n(z) \, dz \right|^2 p(u)p(v) \, du \, dv$$

$$= \frac{1}{h} \int \int \frac{1}{h} \left| M_n\left(\frac{v-u}{h}\right) \right|^2 p(u)p(v) \, du \, dv = T, \qquad \text{say},$$

where $M_n(x) = \int K_n(z+x)K_n(z) \, dz$. Next, use the facts that $p$ is at least $(\beta + \sigma - 1/2)$-Lipschitz continuous and uniformly bounded (Lemma 3 in the Appendix [4]) to find $C$ and $M^Y$, positive constants depending only on $\beta, L$ and $\sigma$, such that

$$\left| T - \frac{1}{h} \|p\|_2^2 \|M_n\|_2^2 \right| \leq \frac{1}{h} \int \int |M_n(x)|^2 |p(v+hx) - p(v)| \, dx \, p(v) \, dv$$

$$\leq \frac{1}{h} \int \left( \int_{|hx| \leq \epsilon} |M_n(x)|^2 C \epsilon^{\beta + \sigma - 1/2} \, dx \right.$$

$$\left. + \int_{|x| > \epsilon/h} 2M^Y |M_n(x)|^2 \, dx \right) p(v) \, dv$$

$$\leq \frac{1}{h} o(\|M_n\|_2^2),$$

where $o(1) \to 0$ as $h \to 0$, depending only on $\beta, L$ and the density $g$. We choose $\epsilon \to 0$ such that $\epsilon/h \to \infty$ so that

$$(15) \qquad T = \frac{\|p\|_2^2 \|M_n\|_2^2}{h} (1 + o(1)).$$

By the Plancherel formula, $\|M_n\|_2^2 = \int |\Phi^{K_n}(u)\Phi^{K_n}(-u)|^2\,du = (\pi(4\sigma+1)\times h^{4\sigma})^{-1}(1+o(1))$. Note that we should again split the integration domain and evaluate the dominant term in the previous integral.

On the other hand, let us deal now with

$$
(16) \quad \begin{aligned}
E_f[\langle K_{n,h}(\cdot - Y_1), K_h \star f\rangle^2] &= E_f\left[\frac{1}{2\pi}\left|\int e^{iuY_1}\frac{\Phi^K(hu)}{\Phi^g(u)}\overline{\Phi}(u)\,du\right|^2\right] \\
&= U, \qquad \text{say.}
\end{aligned}
$$

In Lemma 4 in the Appendix [4] we prove that $\overline{\Phi}/\Phi^g$ is absolutely integrable if $\beta \geq \sigma$. Then by the Lebesgue convergence theorem we see that there exists a function

$$
F(y) = \frac{1}{2\pi}\int e^{iuy}\frac{\overline{\Phi}(u)}{\Phi^g(u)}\,du = \lim_{h\to 0}\frac{1}{2\pi}\int_{|u|\leq 1/h} e^{iuy}\frac{\overline{\Phi}(u)}{\Phi^g(u)}\,du,
$$

which is uniformly continuous and bounded such that $\|F\|_2^2 = \|\overline{\Phi}/\Phi^g\|_2^2/(2\pi)$. Note that $\overline{\Phi}(u) = \Phi(-u)$ and $\Phi^g(u) = \overline{\Phi^g}(-u)$, giving $\overline{F} = F$. Thus, $F$ is a real-valued function.

Thus, we obtain

$$
(17) \qquad\qquad U = E_f[F^2(Y)](1+o(1)).
$$

Finally,

$$
(18) \qquad\qquad \|K_h \star f\|_2^4 = \|f\|_2^4(1+o(1))
$$

by the bias computations.

Thus, from (15), (17) and (18) we get

$$
(19) \quad \begin{aligned}
E_f[S_1^2] &= \frac{2\|p\|_2^2}{\pi(4\sigma+1)}\frac{1+o(1)}{n^2h^{4\sigma+1}} \\
&\quad - \frac{4E_f[F^2(Y)](1+o(1))}{n^2} + \frac{2\|f\|_2^4(1+o(1))}{n^2} \\
&= \frac{2\|p\|_2^2}{\pi(4\sigma+1)}\frac{1+o(1)}{n^2h^{4\sigma+1}}.
\end{aligned}
$$

Use (17) and (18) to get

$$
(20) \quad E_f[S_2^2] = \frac{4}{n}(E_f[F^2(Y)] - \|f\|_2^4)(1+o(1)) = \frac{4\Omega_g^2(f)}{n}(1+o(1)).
$$

The upper bound for the variance follows from (19) and (20) for the case $\beta \geq \sigma$.

For the case $\beta < \sigma$, go back to (16):

$$E_f\left[\frac{1}{2\pi}\left|\int e^{iuY_1}\frac{\Phi^K(hu)}{\Phi^g(u)}\overline{\Phi}(u)\,du\right|^2\right]$$

$$\leq \frac{1}{2\pi}\left(\int_{|u|\leq 1/h}\left|\frac{\overline{\Phi}(u)}{\Phi^g(u)}\right|du\right)^2$$

(21)
$$\leq \frac{1}{2\pi}\left(O(1)+\int_{M\leq|u|\leq 1/h}|u|^\sigma|\Phi(u)|\,du\right)^2$$

$$\leq O(1)\int_{M\leq|u|\leq 1/h}|u|^{2\sigma-2\beta}\,du\cdot\int_{M\leq|u|\leq 1/h}|u|^{2\beta}|\Phi(u)|^2\,du$$

$$\leq \frac{O(1)}{h^{2(\sigma-\beta)+1}}.$$

So, from (19) and (21) we obtain the upper bound for the variance when $\beta < \sigma$. □

An easy consequence of Proposition 1 is that if the underlying unknown density is smoother enough than the noise ($\beta > \sigma + 1/4$) our parameter can be estimated at the parametric rate. We establish next asymptotic normality and a Cramér–Rao type of asymptotic efficiency bound.

THEOREM 1.   *If $\beta > \sigma + 1/4$, the estimator $d_n$ defined in (11) with bandwidth $h = h_*$ such that*

$$n^{-1/(4\sigma+1)} \ll h_* \ll n^{-1/(4\beta)}$$

*is an asymptotically normally distributed estimator of $d$, that is,*

$$\sqrt{n}(d_n - d) \xrightarrow{d} N(0, 4\Omega_g^2(f)).$$

*Moreover, it is asymptotically efficient, attaining the Cramér–Rao bound.*

PROOF.   Let us decompose the risk of the estimator as

$$E_f[|d_n - d|] \leq B(d_n) + \sqrt{V(d_n)} \leq Lh^{2\beta} + \frac{2\Omega_g(f)}{\sqrt{n}}(1 + o(1)),$$

and then use Proposition 1. Indeed, if $\beta > \sigma + 1/4$ and if $n^{-1}h^{-(4\sigma+1)} \ll 1$ we see that $4\Omega_g^2(f)/n(1+o(1))$ is the dominant term in the variance. Let us take $h = o(n^{-1/(4\beta)})$ such that the bias is infinitely smaller, $Lh^{2\beta} \ll 2\Omega_g(f)/n$. So

$$\sqrt{n}(d_n - d) = \sqrt{n}(d_n - E_f[d_n]) + \sqrt{n}B(d_n).$$

The second term of the sum on the right-hand side term tends to 0 and the asymptotic normality of the first term can be deduced from Butucea [3]. It is in this case a classical central limit theorem for $U$-statistics of order 1.

For the Cramér–Rao bound, we follow the lines of proof in Laurent [20]. Similar results were given by Bickel and Ritov [1] following the theory of Ibragimov and Khas'minskii [17] and Levit [23]. A first step of the proof is to compute the Fréchet derivative of the functional $\int f^2 = \int F \cdot p$ at the likelihood $p_0 = f_0 \star g$,

$$\int F \cdot p - \int F_0 \cdot p_0 = \int 2F_0(p - p_0) + \int (F - F_0)(p - p_0)$$

and $\int (F - F_0)(p - p_0) = o(\|p - p_0\|_2)$, when $\|p - p_0\|_2 \to 0$. Next, consider the space orthogonal to the square root of the likelihood $\sqrt{p_0}$, $H = \{k : \int k\sqrt{p_0} = 0\}$ and the projection operator onto this space: $P_{H(p_0)}(k) = k - (\int k\sqrt{p_0})\sqrt{p_0}$. Write $K_n = K = T'(p_0)\sqrt{p_0} = P_{H(p_0)}(k)$ as $\langle g, k \rangle$. Then the minimal variance is $\|g\|_2^2$.

Here, $T'(p_0)k = \int 2F_0 k$; then

$$K = \int 2F_0\sqrt{p_0}\Big(k - \Big(\int k\sqrt{p_0}\Big)\sqrt{p_0}\Big)$$

$$= \int (2F_0\sqrt{p_0})k - \Big(\int 2F_0 p_0\Big)\int \sqrt{p_0}k.$$

So, finally,

$$\|g\|_2^2 = 4\int |F_0|^2 p_0 - \Big(\Big|\int 2F_0 p_0\Big|\Big)^2 = 4V_{f_0}(F_0(Y)). \qquad \square$$

In the following theorem we compute the rate on the nonparametric side ($0 < \beta \le \sigma + 1/4$). We prove in Section 4 that this rate is optimal in the minimax approach under the following additional assumption on the noise distribution.

ASSUMPTION (P).   The distribution of the polynomial noise in (4) is such that $\Phi^g$ is at least three times continuously differentiable. Moreover there exist $A_1$, $A_2 > 1$, $u_0$, $u_1$, $u_2 > 0$ large enough such that

$$|\Phi^g(u)| \ge u_0 \qquad \forall |u| \le A_1$$

and

$$|(\Phi^g)^{(k)}(u)| \le \frac{u_k}{|u|^{\sigma + k}} \qquad \text{for } k = 1, 2, \ \forall |u| \ge A_2.$$

THEOREM 2.   *If $0 < \beta \le \sigma + 1/4$, the estimator $d_n$ of $d$ defined in (11) with bandwidth $h_*$ satisfies the upper bound (8) for the rate $\varphi_n$, where*

$$h_* = n^{-2/(4\beta + 4\sigma + 1)}, \qquad \varphi_n = n^{-4\beta/(4\beta + 4\sigma + 1)}.$$

*Moreover, under Assumption (P) this rate is minimax.*

PROOF OF (8) FOR THEOREM 2. If $0 < \beta \le \sigma + 1/4$, $\|p\|_2^2/(\pi(4\sigma + 1)n^2h_*^{4\sigma+1})$ is the dominant term in the variance, whether $\beta \ge \sigma$ or $\beta < \sigma$. The bandwidth $h_*$ minimizes the bias plus the variance. The upper bound of the normalized mean error is less than $C = \max\{L, \sqrt{2M^p/(\pi(4\sigma + 1))}\}$; see Lemma 3 in the Appendix [4]. □

3.2. *The other setups.* In the case where the densities are smoother than the noise, we can always define the function $F$ as the inverse Fourier transform of $\Phi/\overline{\Phi^g}$. The next theorem gives us the bandwidth $h_*$ so that $d_n$ is an asymptotically normal and efficient estimator.

THEOREM 3. *The estimator $d_n$ defined in* (11) *with bandwidth $h_*$ such that*

$$h_* \ll \left(\frac{\log n}{4\alpha}\right)^{-1/r}$$

*is asymptotically normally distributed and it is asymptotically efficient, attaining the Cramér–Rao bound $4\Omega_g^2$ (see Definition 3):*

(1) *if $f$ belongs to $S(\alpha, r, L)$ and the noise is $\sigma$-polynomially smooth;*

(2) *if $f$ belongs to $S(\alpha, r, L)$ and the noise is exponentially smooth with $r > s$ or with $r = s$ and $\alpha > \gamma$.*

In the case where the noise is exponentially smooth and smoother than the underlying density estimation is always difficult, that is, only nonparametric slower rates are attained. We prove the lower bounds (9), under the following additional assumption, which is not very restrictive.

ASSUMPTION (E). The exponential noise distribution in (5) has a continuously differentiable Fourier transform such that

$$|(\Phi^g)'(u)| \le O(1)|u|^{\mathcal{A}}\exp(-\gamma|u|^s),$$

for large enough $|u|$ and some fixed constant $\mathcal{A} \in R$.

THEOREM 4. *Let the noise be exponentially smooth. The estimator $d_n$ of $d$ defined in* (11) *with bandwidth $h_*$ satisfies the upper bound* (8) *for the rate $\varphi_n$, where:*

(1) *If $f$ belongs to $W(\beta, L)$,*

$$h_* = \left(\frac{\log n}{2\gamma} - \frac{2\beta + 1}{2\gamma s}\log\frac{\log n}{2\gamma}\right)^{-1/s}, \qquad \varphi_n = L\left(\frac{\log n}{2\gamma}\right)^{-2\beta/s};$$

*moreover, under Assumption* (E) *this rate is minimax.*

TABLE 1
*Rates for estimating $d = \int f^2$ from indirect observations*

| $f \backslash g$ | Polynomial: $|u|^{-\sigma}$ | | Exponential: $\exp(-\gamma|u|^s)$ | |
|---|---|---|---|---|
| $W(\beta, L)$ | $\beta < \sigma + 1/4$ | $O(1)n^{-4\beta/(4\beta+4\sigma+1)}$ | $O(1)(\log n/(2\gamma))^{-2\beta/s}$ | |
| | $\beta \geq \sigma + 1/4$ | $2\Omega_g n^{-1/2}$ | | |
| $S(\alpha, r, L)$ | $2\Omega_g n^{-1/2}$ | | $r < s$ | $O(1)\exp(-2\alpha/h_*^r)$ |
| | | | $(r > s)$ | |
| | | | or | |
| | | | $(r = s, \alpha > \gamma)$ | $2\Omega_g n^{-1/2}$ |

NOTE. $h_*$ is the solution of (22).

(2) *If $f$ belongs to $S(\alpha, r, L)$ and, either $r < s$ or $r = s$ and $\alpha \leq \gamma$, $h_*$ is the solution of*

$$(22) \qquad \frac{2\alpha}{h_*^r} + \frac{2\gamma}{h_*^s} = \log n - (\log\log n)^2$$

*and* $\quad \varphi_n = L\exp(-2\alpha/h_*^r)$; *moreover, under Assumption* (E) *this rate is minimax when $r < s$.*

Note that when the density and the noise are both exponentially smooth, the rates are faster than any logarithm but slower than any polynomial in $n$; except when $r = s$ and $\alpha = \gamma$, the rate is nearly parametric, $\varphi_n = c_3(\log n)^{r/2}/\sqrt{n}$ for $h$ the solution of $h^{r-1}\exp(4\alpha/h^r) = cn$.

**4. Goodness-of-fit tests.** Let us give here the convergence rates for the testing procedure in (12) and optimal choice of tuning parameters. The rates are given in Table 2. Note that for setups where we prove the lower bounds for the testing rate we need to assume that the density $f_0$ in the null hypothesis is such that

$$(23) \qquad f_0(x) \geq \frac{c_0}{1 + |x|^2} \qquad \forall x \in \mathbb{R}.$$

Let us note immediately that we have a similar property for $p_0 = f_0 * g$. Indeed, let $A > 1$ be large enough such that $\int_{-A}^{A} g(x)\,dx > 1/2$. Then there is a $c_0^Y > 0$ such that

$$(24) \qquad p_0(x) \geq \int_{-A}^{A} f_0(x - y)g(y)\,dy \geq c_0^Y \min\left\{\frac{1}{A^2}, \frac{1}{|x|^2}\right\} \qquad \forall x \in \mathbb{R}.$$

We choose to work under the assumption (23) for simplicity. Notice that we can as well solve the problem if $f_0$ decays asymptotically like a polynomial (faster than $1/|x|^2$), but for technical reasons we would need to assume

Table 2
*Rates for testing in $\mathbb{L}_2$-norm from indirect observations*

| $f \backslash g$ | Polynomial: $\|u\|^{-\sigma}$ | Exponential: $\exp(-\gamma\|u\|^s)$ | |
|---|---|---|---|
| $W(\beta, L)$ | $O(1)n^{-2\beta/(4\beta+4\sigma+1)}$ | $\sqrt{L}(\log n/(2\gamma))^{-\beta/s}$ | |
| $S(\alpha, r, L)$ | $O(1)(\log n)^{(\sigma+1/4)/r}n^{-1/2}$ | $r < s$ | $\sqrt{L}\exp(-\alpha/h_*^r)$ |
| | | $r > s$ | $O(1)\dfrac{h_*^{(s-1)_-/4}}{\sqrt{n}}\exp(\frac{\gamma}{h_*^s})$ |

NOTE. $h_*$ is defined in (22).

that the characteristic function of the noise is smoother than $\mathcal{C}^1$. Another way of proving the lower bounds consists of assuming (24), which is less restrictive, but then we have to modify the construction of perturbation functions according to the actual asymptotic behavior of $f_0$.

4.1. *Sobolev densities and polynomial noise.* Though two rates were attainable in the same setup for estimating $d$, only one minimax rate for testing is possible. This phenomenon is similar to the case of testing with direct observations.

THEOREM 5. *The test procedure $\Delta_n^*$ defined in (12) for the threshold $t_n$ attains the rate $\psi_n$ and, under Assumption (P) and (23), $\psi_n$ is a minimax rate of testing over the class $W(\beta, L)$, where*

$$h = h_* = n^{-2/(4\beta+4\sigma+1)}, \qquad t_n = \psi_n = n^{-2\beta/(4\beta+4\sigma+1)}.$$

PROOF OF (6) FOR THEOREM 5. Let us bound from above successively the first and second type errors. Note that, for a fixed density $f_0 \in W(\beta, L)$,

$$E_{f_0}[T_n^*] = \|K_h \star f_0 - f_0\|_2^2 = Lh^{2\beta}o(1),$$

similarly to the proof of Proposition 1. In order to compute the variance let us write

$$T_n^* - E_{f_0}[T_n^*] = \frac{1}{n(n-1)}\sum_{k \neq j}\langle K_{n,h}(\cdot - Y_k) - K_h \star f_0, K_{n,h}(\cdot - Y_j) - K_h \star f_0\rangle$$

$$+ \frac{2}{n}\sum_{k=1}^n \langle K_{n,h}(\cdot - Y_k) - K_h \star f_0, K_h \star f_0 - f_0\rangle.$$

Note that the previous sum is null, since for all $k = 1, \ldots, n$,

$$\langle K_{n,h}(\cdot - Y_k) - K_h \star f_0, K_h \star f_0\rangle$$

$$= \frac{1}{2\pi}\int(e^{iuY_k}/\Phi^g(hu) - \Phi_0(u))\Phi^K(hu)\Phi_0(u)\,du$$

$$= \langle K_{n,h}(\cdot - Y_k) - K_h \star f_0, f_0\rangle.$$

Finally, $V_{f_0}[T_n^*] = S\|p_0\|_2^2 n^{-2} h^{-(4\sigma+1)}(1+o(1))$, where $S = 2/(\pi(4\sigma+1))$. So the first type error can be written as

$$P_{f_0}[|T_n^*| \geq \mathcal{C}^* t_n^2] \leq \frac{E_{f_0}[T_n^{*2}]}{\mathcal{C}^{*2} t_n^4} \leq O(1) \frac{S\|p_0\|_2^2}{\mathcal{C}^{*2}} \leq \frac{\xi}{2}$$

for $\mathcal{C}^*$ large enough. For the second type error, consider a density $f$ in $H_1(\mathcal{C}, \psi_n)$. Then $E_f[T_n^*] = \|K_h \star f - f_0\|_2^2$. The bias can be bounded from above as

$$
\begin{aligned}
B[T_n^*] &= |\|K_h \star f - f_0\|_2^2 - \|f - f_0\|_2^2| \\
&= |\|K_h \star f\|_2^2 - \|f\|_2^2 - 2\langle K_h \star f - f, f_0 \rangle| \\
&\leq \frac{1}{2\pi} \int_{|u|>1/h} |\Phi(u)|^2 \, du + \frac{2}{2\pi} \int_{|u|>1/h} |\Phi(u)| \cdot |\Phi_0(u)| \, du \\
&\leq L h^{2\beta}(1+o(1)),
\end{aligned}
$$

since $\int_{|u|>1/h} |u|^{2\beta} |\Phi_0(u)|^2 \, du = o(1)$, for the fixed density $f_0$. In order to evaluate the variance, let us write

$$
\begin{aligned}
T_n^* - E_f[T_n^*] &= \frac{1}{n(n-1)} \sum_{k \neq j} \langle K_{n,h}(\cdot - Y_k) - K_h \star f, K_{n,h}(\cdot - Y_j) - K_h \star f \rangle \\
&\quad + \frac{2}{n} \sum_{k=1}^n \langle K_{n,h}(\cdot - Y_k) - K_h \star f, f - f_0 \rangle \\
&= S_1(f) + S_2(f - f_0),
\end{aligned}
$$

say. As in Proposition 1, the last two terms are uncorrelated, so $V_f[T_n^*] = E_f[|S_1(f)|^2] + E_f[|S_2(f - f_0)|^2]$. Similar computation leads for $h = h_*$ to the upper bound

$$V_f[T_n^*] \leq \frac{SM^p(1+o(1))}{n^2 h^{4\sigma+1}} + \frac{4\Omega_g^2(f - f_0)}{n} I(\beta > \sigma),$$

where $\Omega_g(f - f_0) = \int (F - F_0)^2 p - (\int (f - f_0)f)^2$ and we have used constant $M^p > 0$, such that $\sup_f \|p\|_\infty \leq M^p$, introduced in Lemma 3 in the Appendix [4].

Let us note that whenever $\beta > \sigma$, we find $M > 0$ large enough such that

$$
\begin{aligned}
\Omega_g^2(f - f_0) &\leq \int (F - F_0)^2 p \leq M^p \|F - F_0\|_2^2 \leq \frac{M^p}{4\pi^2} \int \left| \frac{\Phi(u) - \Phi_0(u)}{\Phi^g(u)} \right|^2 du \\
&\leq \int_{|u| \leq M} c_1 M^{2\sigma} |\Phi(u) - \Phi_0(u)|^2 \, du \\
&\quad + \int_{|u|>M} c_2 |u|^{2\sigma} |\Phi(u) - \Phi_0(u)|^2 \, du
\end{aligned}
$$

$$\leq c_3 M^{2\sigma} \|f - f_0\|_2^2 + \frac{c_2}{M^{2(\beta-\sigma)}} \int_{|u|>M} |u|^{2\beta} |\Phi(u)|^2 \, du$$

$$\leq C\|f - f_0\|_2^{2-2\sigma/\beta},$$

where $C$ is a constant depending only on $\beta, L$ and the fixed noise probability density $g$. This inequality is useful for the limit cases in $H_1$ where $\|f - f_0\|_2 \to 0$. So, the second type error can be bounded as

$$P_f[|T_n^*| < \mathcal{C}^* t_n^2] \leq P_f[|T_n^* - E_f[T_n^*]| > \|f - f_0\|_2^2 - \mathcal{C}^* t_n^2 - B[T_n^*]]$$

$$\leq P_f\Big[\frac{|T_n^* - E_f[T_n^*]|}{\sqrt{V_f[T_n^*]}}$$

$$\geq \frac{\|f - f_0\|_2^2 - \mathcal{C}^*\psi_n^2 - Lh^{2\beta}}{c_1(nh^{2\sigma+1/2})^{-1} + c_2\|f - f_0\|_2^{1-\sigma/\beta} n^{-1/2} I(\beta > \sigma)}\Big].$$

Either $0 < \beta \leq \sigma + 1/4$, when the probability above is less than $c_1^2(\mathcal{C} - \mathcal{C}^* - L)^{-2} \leq \xi/2$ for $\mathcal{C} > \mathcal{C}^*$ large enough, or $\beta > \sigma + 1/4$, when

$$P_f[|T_n^*| < \mathcal{C}^* t_n^2] \leq P_f\Big[\frac{|T_n^* - E_f[T_n^*]|}{\sqrt{V_f[T_n^*]}} \geq c_2\sqrt{n}\psi_n^{1+\sigma/\beta}\Big]$$

$$\leq n^{-(2\beta+2\sigma+1)/(4\beta+4\sigma+1)} = o(1)$$

for $\mathcal{C} > \mathcal{C}^*$ large enough.

The upper bounds in (6) are proved. For the lower bounds in (7) see Section 5. $\square$

### 4.2. *The other setups.*

We know now that in some setups we can estimate $d$ at the parametric $n^{-1/2}$ rate. We shall see next that the minimax testing rate is necessarily nonparametric.

THEOREM 6. *The test procedure $\Delta_n^*$ defined in* (12) *for the bandwidth $h_*$, the threshold $t_n$ and the constant $\mathcal{C}^*$ satisfies the upper bound* (6) *for the rate $\psi_n$, where:*

(1) *If $f$ belongs to $S(\alpha, r, L)$ and the noise is polynomially smooth,*

$$h = h_* = \Big(\frac{\log n}{2\alpha} - \frac{2\sigma + 1/2}{2\alpha r}\log\log n\Big)^{-1/r},$$

$$t_n = \psi_n = \frac{1}{\sqrt{n}}\Big(\frac{\log n}{2\alpha}\Big)^{(4\sigma+1)/(4r)}.$$

(2) *If $f$ belongs to $W(\beta, L)$ and the noise is exponentially smooth,*

$$h = h_* = \left( \frac{\log n}{2\gamma} - \frac{2\beta + 1}{2\gamma s} \log \frac{\log n}{2\gamma} \right)^{-1/s},$$

$$t_n = \psi_n = \sqrt{L} \left( \frac{\log n}{2\gamma} \right)^{-\beta/s};$$

*moreover, under Assumption* (E), *$\psi_n$ is an exact minimax rate of testing.*

(3) *If $f$ belongs to $S(\alpha, r, L)$ and the noise is exponentially smooth, $h = h_*$ is a solution of* (22) *and*

$$t_n = \psi_n = \begin{cases} \sqrt{L} \exp\left( -\dfrac{\alpha}{h_*^r} \right), & \text{if } r < s, \\[2ex] \dfrac{h_*^{(s-1)_-/4}}{\sqrt{n}} \exp\left( \dfrac{\gamma}{h_*^s} \right), & \text{if } r \geq s; \end{cases}$$

*moreover, under Assumption* (E), *$\psi_n$ is an exact ($\mathcal{C}^* = \mathcal{C}_* = 1$) minimax rate of testing for the case $r < s$.*

We prove in the Appendix [4] exact lower bounds for the case $r < s$, but the same proof provides lower bounds precisely within a logarithmic factor for the case $r > s$.

**5. Lower bounds.** We show in the first part that proofs for minimax lower bounds for the estimation problem of $d$ and for the testing problem in $\mathbb{L}_2$ come down to the same choice of hypotheses and to checking similar conditions.

LEMMA 1. *Let $f_0$ and $f_1$ be two probability densities in the class $W(\beta, L)$, depending on $n$. If:*

(a) *for estimation densities are such that $|\|f_1\|_2^2 - \|f_0\|_2^2| \geq 2\varphi_n$, for some $\varphi_n > 0$;*

(a′) *for testing densities are such that $\|f_1 - f_0\|_2 \geq \mathcal{C}\psi_n$, for some $\psi_n > 0$;*

(b) *$P_1^Y \ll P_0^Y$ and there exists $0 < \eta < 1$ such that*

$$\chi^2(P_0^Y, P_1^Y) \stackrel{\text{def}}{:=} \int \left( \frac{dP_1^Y}{dP_0^Y} - 1 \right)^2 dP_0^Y \leq \eta^2,$$

*then*

$$\inf_{\hat{d}_n} \sup_{f \in W(\beta, L)} \varphi_n^{-1} E_f[|\hat{d}_n - d|] \geq (1 - \eta)(1 - \sqrt{\eta}),$$

$$\inf_{\Delta_n} \sup_{f \in W(\beta, L)} \left( P_{H_0}(\Delta_n = 1) + P_{H_1(\mathcal{C}, \psi_n)}(\Delta_n = 0) \right) \geq (1 - \eta)(1 - \sqrt{\eta}).$$

Proof.    For the estimation problem the risk $\sup_{f \in W(\beta,L)} \varphi_n^{-1} E_f[|\hat{d}_n - d|]$ is bounded from below by the risk for two hypotheses, $\max_{i=0,1} \varphi_n^{-1} E_{f_i}[|\hat{d}_n - d_i|]$, and then we directly use Lemma 4 from Butucea and Tsybakov [6].

For the testing problem, we choose two hypotheses, $f_0$, the density under $H_0$, and another density $f_1$ under $H_1$ (which implies that $\|f_1 - f_0\|_2 \geq \mathcal{C}\psi_n$, for some $\psi_n > 0$). Then the risk for the test problem becomes

$$Rtest := \inf_{\Delta_n} \sup_{f \in W(\beta,L)} \left(P_{H_0}(\Delta_n = 1) + P_{H_1(\mathcal{C},\psi_n)}(\Delta_n = 0)\right)$$

$$\geq \inf_{\Delta_n} \left(P_{f_0}^Y(\Delta_n = 1) + (1 - \sqrt{\eta})P_{f_0}^Y\left(\Delta_n = 0, \frac{dP_{f_1}^Y}{dP_{f_0}^Y} \geq 1 - \sqrt{\eta}\right)\right).$$

This gives

$$Rtest \geq (1 - \sqrt{\eta})P_{f_0}^Y\left(\frac{dP_{f_1}^Y}{dP_{f_0}^Y} \geq 1 - \sqrt{\eta}\right)$$

$$\geq (1 - \sqrt{\eta})\left(1 - \frac{1}{\eta}E_{f_0}\left[\left(\frac{dP_{f_1}^Y}{dP_{f_0}^Y} - 1\right)^2\right]\right),$$

which allows one to conclude when assumption (b) holds.   □

We shall use in the proofs the following construction. Let $0 < \delta < 1$ be small through the remaining proofs of lower bounds.

In the estimation problem, let us choose $f_0$, a density function in the Sobolev class $W(\beta, a(\delta)L)$, respectively, $S(\alpha, r, a(\delta)L)$, where $0 < a(\delta) < 1$ is a constant depending on $\delta$ and defined for each setup, such that (23) holds. Moreover, for the estimation problem we want to choose the Fourier transform $\Phi_0$ to have compact support included in $(-2\delta, 2\delta)$.

In the testing problem, we have to assume that the density $f_0$ satisfies (23).

Proof of (9) in Theorem 2 and of (7) in Theorem 5.    This proof is based on a large family of hypotheses. Similar reasoning proves that the same construction is valid for proving lower bounds for both quadratic functional estimation and nonparametric testing in $\mathbb{L}_2$.

Note that this setup includes Theorem 2 for $\beta < \sigma + 1/4$. This is not a contradiction, since the lower bounds here are much slower than the parametric $n^{-1/2}$ rate that the estimator attains; see Theorem 1.

Let $\theta_j, j = 1, \ldots, M$, be independent Bernoulli random variables and let $\Pi$ be the probability measure associated with them. For $h > 0$ small as $n \to \infty$ and for a function $H$ to be defined later, let

$$(25) \qquad f_\theta(x) = f_0(x) + \sum_{j=1}^{M} \theta_j h^{\beta+\sigma+1} H_h(x - x_j),$$

where $H_h(\cdot) = 1/h H(\cdot/h)$, $x_j = jh$ and $M$ is an integer such that $M/h = 1 - o(1)$, as $n \to \infty$ and for $h$ small. Note that the observations $Y_i$, $i = 1, \dots, n$, when the underlying density is $f_\theta$, have density $p_\theta(x) = p_0(x) + \sum_{j=1}^M \theta_j h^{\beta+\sigma+1} G_h(x - x_j)$, where the function $G$ is defined in Lemma 2 and $H$ is such that $\Phi^G(u) = \Phi^H(u)\Phi^g(u/h)$. Indeed, $(H_h(\cdot - x_j) * g)(x) = H_h * g(x - x_j) = G_h(x - x_j)$. Using Lemmas 5 and 6 in the Appendix [4] we see that the hypotheses fit into the model, that is, $f_\theta$ are density functions for all $\theta$ belonging to the Sobolev class $W(\beta, L)$ and such that

$$\Pi[\|f_\theta - f_0\|_2^2 \geq \mathcal{C} n^{-4\beta/(4\beta+4\sigma+1)}] \to 1,$$

as $n \to \infty$, for fixed $\mathcal{C} > 0$.

LEMMA 2.   *Let the function* $G : [-1, 0] \to \mathbb{R}$ *be defined by*

$$G(x) = \exp\left(-\frac{1}{1 - (4x+3)^2}\right) I_{[-1, -1/2]}(x)$$
$$- \exp\left(-\frac{1}{1 - (4x+1)^2}\right) I_{[-1/2, 0]}(x).$$

*Then* $G$ *is an infinitely differentiable function such that* $\int G(x)\, dx = 0$ *and having all polynomial moments finite. Its Fourier transform is such that*

$$|\Phi^G(u)| \leq C_G \exp(-a\sqrt{|u|}) \qquad as\ |u| \to \infty$$

*for some positive constants* $C_G$, $a > 0$. *Moreover* $\Phi^G$ *is an infinitely differentiable, bounded function.*

This construction is based on the function $f_a$ in Lepski and Levit [21], page 133, and the asymptotic behavior of its Fourier transform follows from the reference therein.

We stress the fact that in this setup hypotheses functions $f_\theta$ belong to $H_1(\mathcal{C}, \psi_n)$ with probability which tends to 1 when $n \to \infty$. In order to bound the risk from below, very small modification is needed in the proof of Lemma 1 that we do not discuss in detail here. The last thing to check is that the distance between resulting models is finite,

$$\Delta^2 := E_{f_0}\left[\left(\frac{\int \prod_{i=1}^n p_\theta(Y_i)\pi(d\theta) - \prod_{i=1}^n p_0(Y_i)}{\prod_{i=1}^n p_0(Y_i)}\right)^2\right]$$
$$= E_{f_0}\left[\left(\int \prod_{i=1}^n \left(1 + \sum_{j=1}^M \theta_j h^{\beta+\sigma+1} \frac{G_h(Y_i - x_j)}{p_0(Y_i)}\right)\pi(d\theta_j)\right)^2\right] - 1.$$

Now, denote by $Y_{i,j}$ those observations $Y_i$ belonging to the support of $G_h(\cdot - x_j)$ and $a_{ij} = h^{\beta+\sigma+1} G_h(Y_{i,j} - x_j)/p_0(Y_{i,j})$. Since those intervals are disjoint

we write

$$\Delta^2 = E_{f_0}\left[\left(\int \prod_{i=1}^{n}\prod_{j=1}^{M}\left(1 + \theta_j h^{\beta+\sigma+1}\frac{G_h(Y_{i,j}-x_j)}{p_0(Y_{i,j})}\right)\pi(d\theta_j)\right)^2\right] - 1$$

$$= E_{f_0}\prod_{j=1}^{M}\left\{\frac{1}{2}\prod_{i=1}^{n}(1+a_{ij}) + \frac{1}{2}\prod_{i=1}^{n}(1-a_{ij})\right\}^2 - 1$$

$$\leq E_{f_0}\prod_{j=1}^{M}\left\{\frac{1}{2}\prod_{i=1}^{n}(1+a_{ij}^2) + \frac{1}{2}\prod_{i=1}^{n}(1-a_{ij}^2)\right\} - 1,$$

where we have used the facts that $(a+b)^2 \leq 2a^2 + 2b^2$ and that $E_{f_0}[a_{ij}] = 0$ since $\int G = 0$. Moreover, $a_{ij}$ are small with $n$ and $E_{f_0}[a_{i,j}^2] \leq ch^{2\beta+2\sigma+1}$ by Lemma 6 in the Appendix [4]. Therefore

$$\Delta^2 \leq E_{f_0}\left[\prod_{j=1}^{M}\left(1 + \sum_{i_1 \neq i_2} a_{i_1 j}^2 a_{i_2 j}^2\right)\right] - 1$$

$$\leq \sum_{j=1}^{M} n(n-1)E_{f_0}[a_{i_1 j}^2 a_{i_2 j}^2],$$

which is smaller than $cMn^2 h^{4\beta+4\sigma+2} < c'$. $\quad\square$

**Acknowledgments.** I would like to thank two anonymous referees, an Associate Editor and the Editor for their careful reading and suggestions which considerably improved the manuscript.

LABORATOIRE PAUL PAINLEVÉ, UMR CNRS 8524
UFR DE MATHÉMATIQUES
59655 VILLENEUVE D'ASCQ CEDEX
FRANCE
E-MAIL: cristina.butucea@math.univ-lille1.fr